\documentclass[12pt,twoside]{amsart}

\renewcommand{\c}[0]{{\mathbb C}}

\renewcommand{\r}[0]{{\mathbb R}} 

\renewcommand{\a}[0]{{\mathbb A}}

\newcommand{\p}[0]{{\mathbb P}}

\newcommand{\map}[0]{\dasharrow}
\newcommand{\qtq}[1]{\quad\mbox{#1}\quad}
\newcommand{\spec}[0]{\operatorname{Spec}}

\newcommand{\gal}[0]{\operatorname{Gal}}

\newcommand{\im}[0]{\operatorname{im}}

\newsymbol\subsetneq 2328
\newsymbol\onto 1310

\newsymbol\twoheadrightarrow 1310

\newtheorem{thm}{Theorem}[section]

\newtheorem{lem}[thm]{Lemma}
\newtheorem{cor}[thm]{Corollary}

\newtheorem{prop}[thm]{Proposition}

\theoremstyle{definition}

\newtheorem{say}[thm]{}
\newtheorem{exmp}[thm]{Example}

\newtheorem{rem}[thm]{Remark}          
\newtheorem{ack}{Acknowledgments}        
   
\newtheorem{defn-thm}[thm]{Definition--Theorem}

\theoremstyle{remark}

\begin{document}
\bibliographystyle{amsplain}

\title{Fundamental groups of rationally connected varieties}
\author{J\'anos Koll\'ar}

\maketitle

\section{Introduction}

Let $X$ be a smooth, projective, unirational  variety and 
$U\subset X$ an open set.
 The aim of this paper is to find a smooth rational curve $C\subset X$
such that the fundamental group of $C\cap U$ surjects onto the
fundamental group of $U$. 
Following the methods of \cite{rcloc} and \cite{ct-harb},
  a positive answer
over $\c$ translates  to a positive answer over any $p$-adic field.
This gives a rather geometric proof of the theorem of
\cite{harb} about the existence of Galois covers of the  line over $p$-adic
 fields (\ref{harb.cor}). We also obtain a slight
generalization of  the results of \cite{ct-harb} about the
existence of certain torsors  over open subsets of the  line
over $p$-adic fields (\ref{tors.gen.cor}).

If $U=X$ then $\pi_1(X)$ is trivial  (cf.\ (\ref{rc->simp.conn})),
 thus any rational
curve $C$ will do. If $X\setminus U$ is a divisor with normal crossings
and $C$ intersects every irreducible component of  $X\setminus U$
transversally, then the {\it normal} subgroup of $\pi_1(U)$ generated by 
the image of 
$\pi_1(C\cap U)$ equals $\pi_1(U)$ by a simple argument. 
(See, for instance, the beginning of (\ref{pf.main.thm}).)
It is also not hard to produce rational curves $C$ such that the image of
$\pi_1(C\cap U)$ has finite index in  $\pi_1(U)$  
(cf.\ (\ref{maps.fund.grps})). 
These  results  suggest that we are very close to a complete answer,
but surjectivity is not obvious.
Differences between surjectivity and finiteness of the  index 
appear in many similar situation; see, for instance,
\cite[Part I]{shafbook} or \cite{nap-ram}.

The present proof relies on  the  machinery of rationally connected varieties
developed in the papers \cite{kmm1, kmm2, kmm3}.  
The relevant facts are recalled in Section 2.
\medskip

The main geometric result is the following:

\begin{thm}\label{main.thm}
 Let $K$ be an algebraically closed field of characteristic zero
and  $X$ a smooth, projective variety over $K$
which is  rationally connected  (\ref{rc.equiv.def}).
Let $U\subset X$ be an open subset and $x_0\in U$ a point. 
Then there is an open
 subset $0\in V\subset \a^1$ and a morphism $f:V\to U$ such 
that $f(0)=x_0$ and
$$
\pi_1(V,0)\onto \pi_1(U,x_0)\qtq{is surjective.}
$$
Moreover, we can  assume that the following also hold:
\begin{enumerate}
\item $H^1(\p^1, \bar f^*T_X(-2))=0$ where $\bar f:\p^1\to X$ is the unique
 extension of $f$.
 \item $\bar f$ is an embedding if $\dim X\geq  3$ and an
 immersion if $\dim X=2$.
 \end{enumerate}
\end{thm}

\begin{cor}\label{main.cor}
 Let $K$ be a $p$-adic  field 
and  $X$ a smooth, projective variety over $K$ which is 
rationally connected  over
 $\bar K$.
Let $U\subset X$ be an open subset and $x_0\in U(K)$ a point. Then 
there is an open
 subset $0\in V\subset \a^1$ and a morphism $f:V\to U$ 
 (all defined over $K$) such that $f(0)=x_0$ and
$$
\pi_1 (V,0)\onto \pi_1 (U,x_0)\qtq{is surjective,}
$$
where $\pi_1$ here denotes the algebraic fundamental group.
\end{cor}

\begin{rem}\label{large.field.defn}
 More generally, (\ref{main.cor}) holds for any field $K$
of characteristic zero
 such that every curve with a smooth
$K$-point contains a Zariski dense set of $K$-points.
Characterizations of this property are
given in \cite[1.1]{pop}.
The following
are some interesting classes of such fields:
\begin{enumerate}
\item Fields complete with respect to a discrete valuation.
\item Quotient fields of   local Henselian domains.
\item $\r$ and all real closed fields.
\item Pseudo algebraically closed fields, cf.\ \cite[Chap.\
10]{fj86}
\end{enumerate}
\end{rem}

\begin{cor}\cite{harb}\label{harb.cor}
 Let $G$ be a finite group and $K$  a   field of characteristic zero
 as in (\ref{large.field.defn}). Then there
  is a  Galois cover $g:C\to \p^1_K$ with Galois group $G$
such that $C$ is geometrically irreducible and 
 $g^{-1}(0{\colon\!} 1)\cong G$.
\end{cor}

Proof. Let $G\subset GL(n,K)$ be a faithful representation.
Set $U=GL(n)/G$ with 
quotient map $h:GL(n)\to U$ and let 
$x_0$ be the image of the
identity matrix. Then $U$ is unirational, thus by
(\ref{main.cor}) there is a $0\in V\subset \a^1$
and a morphism $f:V\to U$ such that $\pi_1 (V)\onto\pi_1 (U)$
is onto. $h:GL(n)\to U$ is  \'etale and proper, thus it corresponds to a
quotient  $\pi_1 (U)\onto G$. 
 The fiber product
$W:=GL(n)\times_U V\to V$  corresponds to the surjective homomorphism 
$$
\pi_1 (V)\onto\pi_1 (U)\onto G.
$$
Thus $W$  is connected and 
 $W\to V$ is a Galois cover with Galois group $G$.
 Since $W$ has a $K$-point, it is also geometrically connected.
 The preimage of $0\in V$ is  isomorphic to $G$
 (the disjoint union of $|G|$ copies of $\spec_K$).
$W\to V$ can be extended to a (ramified) 
Galois cover of the whole $\p^1_K$. 
\qed

\begin{rem} The above proof works in positive characteristic
if we know that for every subgroup $H<G$ the quotient
$GL(n)/H$ has a smooth compactification.
\end{rem}

The following result was proved by \cite{ct-harb} for
finite groups, which is probably the most important for applications.

\begin{cor}\label{tors.gen.cor}
 Let $K$ be a field of characteristic zero as in (\ref{large.field.defn}),
$G$  a linear algebraic group scheme over $K$ and  $A$  a
principal homogeneous $G$-space.
 Then there is an open set
$0\in V\subset \a^1_K$ and a geometrically irreducible
$G$-torsor $g:W\to V$ such that $g^{-1}(0)\cong A$
(as a $G$-space).
\end{cor}

Proof. Assume that $G$ acts on $A$ from the left
and choose an embedding $G\subset GL(n)$ over $K$.
 $A\times GL(n)$ has a diagonal left action by $G$
and a right action by $GL(n)$ acting only on $GL(n)$. 

The right $G$-action makes the morphism
 $h:  G\backslash(A\times GL(n))\to 
 G\backslash(A\times GL(n))/G=:U$ into a $G$-torsor.
Let $x_0\in U$ be the image of $G\backslash(A\times G)$. 
The fiber of $h$ over $x_0$ is isomorphic to
$A$. Let $G^0$ be the connected component of $G$. 
Then $G\backslash(A\times GL(n))\to G\backslash(A\times GL(n))/G^0$
is smooth with connected fibers
and $G\backslash(A\times GL(n))/G^0\to U$ is \'etale and proper.
Let $0\in V\subset \a^1$  
and $f:V\to U$ be as in the proof of (\ref{harb.cor}).
Then  $W:=(G\backslash(A\times GL(n)))\times_UV$ works.\qed

\medskip

A geometric application is the following:

\begin{cor} For every $2\leq g\leq 13$ there is an open set
$0\in V_g\subset \c$ and a smooth proper morphism
with genus $g$ fibers $S_g\to V_g$ such that the
image of the  monodromy representation
is the full Teichm\"uller group.
\end{cor}

Proof. The moduli of curves is unirational for $g\leq 13$.
Apply (\ref{main.thm}) to the open subset of
curves without automorphisms  $U_g\subset M_g$.\qed

\section{Rationally connected varieties and morphisms of rational curves}

Rationally connected varieties were introduced in \cite{kmm2}
as a higher dimensional generalization of rational and unirational
varieties.
A surface is rationally connected iff it is rational.
In higher dimensions  rationality
 and unirationality are very hard to check.
The notion of rational connectedness  concentrates on
rational curves on a variety.
The following characterizations
were developed in the papers
\cite{kmm2}, \cite[IV.3]{rcbook}, \cite[4.1.2]{k-lowdeg}.

\begin{defn-thm}
 \label{rc.equiv.def}
Let $ K$ be an  algebraically closed field of characteristic zero.
 A smooth proper variety $X$ over $ K$
is called {\em rationally connected} if it satisfies any of the
following 
equivalent properties:
\begin{enumerate}
\item There is an open subset $\emptyset\neq X^0\subset X$,
such that for every
$x_1,x_2\in X^0$, there is a morphism $f:\p^1\to X$ satisfying
$x_1,x_2\in f(\p^1)$.

\item For every
$x_1,\dots,x_n\in X$, there is a morphism $f:\p^1\to X$
satisfying
$x_1,\dots,x_n\in f(\p^1)$.

\item There is a morphism $f:\p^1\to X$ such that 
$H^1(\p^1,f^*T_X(-2))=0$. (This is equivalent to $f^*T_X$ being 
ample.)

\item There is a variety $P$ and a dominant morphism
$F:\p^1\times P\to X$ such that $F((0{\colon\!}1)\times P)$ is a point.
We can also assume that $H^1(\p^1,F_p^*T_X(-2))=0$ for every $p$ where
$F_p:=F|_{\p^1\times\{p\}}$.

\item  Let $z_1,\dots,z_n\in\p^1$ be distinct points
and $m_1,\dots,m_n$ natural numbers. For
each $i=1,\dots,n$ let
$f_i:\spec  K[t]/(t^{m_i})\to X$ be a morphism.  
 Then there is a variety $P$ and a dominant morphism
$F:\p^1\times P\to X$
 such that 
 \begin{enumerate}
 \item the Taylor series of $F_p$ at $z_i$ coincides
with $f_i$ 
up to order $m_i$ for every $i$ and $p\in P$,
\item $H^1(\p^1,F_p^*T_X(-\sum m_i))=0$ for every $i$ and $p$.\qed
\end{enumerate}
\end{enumerate}
\end{defn-thm}

Another easy result that we need is the following.

\begin{lem}(cf.\ \cite[II.3.5.4,II.3.10.1 and II.3.11]{rcbook})
\label{free.maps.lem}
Let $X, P$ be  smooth varieties and 
$F:\p^1\times P\to X$ a dominant morphism
 such that $F((0{\colon\!} 1)\times P)$ is a point.
 Then there is a dense open set $P^0\subset P$ such that 
 $H^1(\p^1,F_p^*T_X(-2))=0$ for every $p\in P^0$ where
$F_p:=F|_{\p^1\times\{p\}}$.

Conversely, let $f:\p^1\to X$ be a morphism such that 
$H^1(\p^1,f^*T_X(-2))=0$. Then there is a pointed variety $p_0\in P$
and a dominant morphism $F:\p^1\times P\to X$ 
such that $F((0{\colon\!} 1)\times P)$ is a point,
$F$ is smooth away from $(0{\colon\!} 1)\times P$ 
and $F_{p_0}=f$.
\qed
\end{lem}

The following result was proved by 
\cite{serre} for unirational varieties and  by \cite{campana-twist, kmm3}
 in general.

\begin{prop}\label{rc->simp.conn}
A  smooth, proper, rationally connected variety is simply connected.\qed
\end{prop}

In the course of the proof we repeatedly encounter the following
situation.  We have morphisms $f^i:\p^1\to X$
each passing through the same point $x_0\in X$.
We would like to have a family of morphisms $f_t:\p^1\to X$
such that the union of the maps $f^i$ can be considered as the limit
of the maps $f_t$ as $t\mapsto 0$. The following lemma is a
technical formulation of this idea. Its statement is  a bit complicated 
since we also want to keep track of the field over
 which the $f_t$ are defined.

\begin{lem}(cf. \cite[3.2]{rcloc})
\label{smoothing.lem} Let $K$ be a field
and $x_0\in X$ a  smooth, proper, pointed  $K$-scheme. Let
$S$ be a zero dimensional reduced $K$-scheme
and $f_0:\p^1_S\to X$ a morphism such that
\begin{enumerate}
\item  $H^1(\p^1_S,f_0^*T_X(-2))=0$, and
\item $f_0(S\times \{(0{\colon\!} 1)\})=\{x_0\}$.
\end{enumerate}
Then there are
\begin{enumerate}\setcounter{enumi}{2}
\item   a smooth pointed curve $0\in D$ over $K$,
\item  a smooth surface $Y$ with a proper morphism
$h:Y\to D$ and a section $B\subset Y$ of $h$, and
\item a morphism $F:Y\to X$,
\end{enumerate}
such that
\begin{enumerate}\setcounter{enumi}{5}
\item    $h^{-1}(0)$ is the union of $\p^1_S$ with a copy 
$B_0$ of $\p^1_K$ such that  $B_0\cap \p^1_S=
S\times \{(0{\colon\!} 1)\}$ and $B_0\cap B$ is a single point.
\item   $F$ restricted to $\p^1_S$ coincides with $f_0$
and $F(B_0\cup B)=\{x_0\}$,
\item  $h^{-1}(D^0)\cong \p^1_K\times  D^0$, where $D^0:=D\setminus \{0\}$,
\item $H^1(\p^1_d, F_d^*T_X(-2))=0$  for every
$d\in D^0$ where $F_d:=F|_{h^{-1}(d)}$.
\end{enumerate}
\end{lem}

Proof. Let us start with any curve $0\in D'$ and
$Y':=\p^1\times D'$. Let $S'\subset \p^1\times \{0\}$ be a subscheme
isomorphic to $S$ and   $Y''$   the blow up of $S'\subset D'$
with projection $h':Y''\to D'$. We can define a morphism
$f':(h')^{-1}(0)\to X$ by setting $f'$ to be $f_0$
on the exceptional divisor of $Y''\to Y'$ and the constant morphism
 to $\{x_0\}$
on the birational transform of $\p^1\times \{0\}$. 
Fix any section $B'\subset Y'$ which does not pass through $S'$. 
We are done if $f'$ can be extended to $F':Y''\to X$
as required. In general this is not possible, but
 such an extension exists after
a suitable \'etale base change $(0\in D)\to (0\in D')$. 
This is proved in \cite[1.2]{kmm2} and \cite[2.2]{rcloc}.\qed

\section{Fundamental groups of fibers of morphisms}

We need some easy results about the variation of fundamental groups
for fibers of nonproper morphisms.

\begin{lem}\label{open.fibr.lem}
 Let $K$ be an algebraically closed 
 field of characteristic zero, $\bar Z,D$  irreducible
$K$-varieties and $f:\bar Z\to D$ a smooth and proper morphism with
 connected fibers. Let $Z\subset \bar Z$ be an open subset
such that $(\bar Z\setminus Z)\to D$ is smooth. 
Let $z_0\in Z(K)$  be a $K$-point, $d_0=f(z_0)$ and $Z_0$ the fiber of
$Z\to D$ through $z_0$.
Then there is an exact sequence
$$
\pi_1(Z_0,z_0)\to \pi_1(Z,z_0)\to \pi_1(D,d_0)\to 1.
$$
\end{lem}

Proof. Over $\c$ the fibration $Z(\c)\to D(\c)$
is a topological fiber bundle, thus we have the above exact
sequence. To settle the algebraic case, let $Y\to Z$ be any
connected 
finite degree \'etale cover and extend it to a finite morphism
$\bar Y\to \bar Z$ where $\bar Y$ is normal. Since
$\bar Z\setminus Z$ is smooth over $D$, the same holds for
$\bar Y\to D$. 
(This is a special case of Abhyankar's lemma,
cf.\ \cite[XIII.5.2]{sga1}.)
The generic fiber of $\bar Y\to D$
is irreducible. Since $\bar Y\to D$
 is smooth and proper, every fiber is irreducible.
This is equivalent to the exactness of the above sequence.\qed
\medskip

The following technical lemma is an upper semi-continuity
statement for the fundamental groups of fibers of nonproper
morphisms.

\begin{lem} \label{fg.uppsem.lem}
 Let $K$ be an algebraically closed 
 field of characteristic zero, $W$ a normal surface over $K$,
$f:W\to D$ a (not necessarily proper)
morphism to a curve with  connected fibers.
Let $B\subset W$ be  a connected subset, one of whose
 irreducible components is a
 section of $f$. Let $ d_0\in D$ be a $K$-point and $C_0$ an
   irreducible component of $ f^{-1}(d_0)$ with a
$K$-point $b_0\in C_0\cap B$. Let 
$x_0\in U$ be a pointed $K$-scheme and 
$h:W\to U$  a morphism such that
$h(B)=\{x_0\}$. Then there is an open subset $D^0\subset D$
such that for every $d\in D^0$ and for every $K$-point
$b_d\in C_d\cap B$, 
$$
\im[\pi_1(C_d,b_d)\to \pi_1(U,x_0)]\supset
\im[\pi_1(C_0,b_0)\to \pi_1(U,x_0)].
$$
\end{lem} 

Proof. Choose a normal compactification $\bar f:\bar W\to D$.
Let $D^0\subset D$ be any open subset such that
 $\bar f$ is smooth with  irreducible fibers over $D^0$ and
  $\bar W\setminus W\to D$ is unramified over $D^0$.
Set $W^0:=f^{-1}(D^0)$. 
$\pi_1(C_0,b_0)\to \pi_1(U,x_0)$ factors through
$\pi_1(W,b_0)\to \pi_1(U,x_0)$ and it has the same image as
$\pi_1(W,b)\to \pi_1(U,x_0)$ for any $b\in B(K)$. 
$\pi_1(W^0,b)\to\pi_1(W,b)$ is surjective by (\ref{maps.fund.grps}).
 By (\ref{open.fibr.lem})
there is an exact sequence
$$
\pi_1(C_d,b_d)\to \pi_1(W_0,b_d)\to \pi_1(D^0,d)\to 1.
$$
Let $B^0\subset B\cap W^0$ be a section
of $f$. Then $\pi_1(B^0,b_d)$ maps onto
$\pi_1(D^0,d)$ and the image of
$\pi_1(B^0,b_d)$ in $\pi_1(U,x_0)$ is trivial.
Thus 
$$
\im[\pi_1(C_d,b_d)\to \pi_1(U,x_0)]=
\im[\pi_1(W^0,b_d)\to \pi_1(U,x_0)]
$$
and we are done.\qed

\begin{lem}(cf.\ \cite[IX.5.6]{sga1}, \cite[1.3]{campana-twist}).
\label{maps.fund.grps}
Let $X,Y$ be  normal varieties, $x_0\in X$  a closed 
point  and $f:X\to Y$ a dominant morphism.
 Then $\im[\pi_1(X,x_0)\to \pi_1(Y,f(x_0))]$
has finite index in $\pi_1(Y,f(x_0))$. If $f$ is an open immersion
then $\pi_1(X,x_0)\to \pi_1(Y,f(x_0))$ is surjective.\qed
\end{lem}

\section{Proof of the main results}

The theory of free morphisms of curves
(cf.\ \cite[II.3]{rcbook}) suggests that
 morphisms
$f:\p^1\to X$ such that $H^1(\p^1, f^*T_X(-2))=0$ behave rather predictably.
Therefore we concentrate on such morphisms. First we establish that
there is a unique maximal subgroup of $\pi_1(U)$ obtainable 
from such a morphism.

\begin{lem} \label{biggest.sg.lem}
Let $X$ be a smooth, proper, rationally connected variety
over an algebraically closed field of characteristic zero. Let
$U\subset X$ be an open set and $x_0\in U$ a point. Then there
is a unique finite index subgroup $H<  \pi_1(U,x_0)$
with the following properties:
\begin{enumerate}
\item There is  a morphism
$f:\p^1\to X$ such that $H^1(\p^1, f^*T_X(-2))=0$,
 $f(0{\colon\!} 1)=x_0$  and
$H=\im[\pi_1(f^{-1}(U),(0{\colon\!} 1))\to \pi_1(U,x_0)]$. 
\item Let $g$ be any morphism
$g:\p^1\to X$ such that $g(0{\colon\!} 1)=x_0$ and 
 $H^1(\p^1, g^*T_X(-2))=0$. Then
$H\supset \im[\pi_1(g^{-1}(U),(0{\colon\!} 1))\to \pi_1(U,x_0)]$. 
\end{enumerate}
\end{lem}

Proof.  First we find one morphism as in 
(\ref{biggest.sg.lem}.2) such that 
$$
\im[\pi_1(g^{-1}(U),(0{\colon\!} 1))\to \pi_1(U,x_0)]
$$
 has finite index in 
$\pi_1(U,x_0)$.

Let $F:\p^1\times P\to X$ be as in 
 (\ref{rc.equiv.def}.4). 
Let $P^0\subset P$ be an open subset such that
$F^{-1}(X\setminus U)\to P$ is \'etale over $P^0$. 
Pick any poit $p\in P^0$. 
By
(\ref{open.fibr.lem}) there is an exact sequence
$$
\begin{array}{lll}
\pi_1(F_p^{-1}(U), (0{\colon\!} 1)\times \{p\}))
&\to & \pi_1((\p^1\times P^0)\cap F^{-1}(U), (0{\colon\!} 1)\times \{p\}))\\
& \to &\pi_1(P^0,p)\to 1.
\end{array}
$$
The section $(0{\colon\!} 1)\times P^0$ is mapped to  a point by $F$,
thus
$$
\pi_1(F_p^{-1}(U), (0{\colon\!} 1)\times \{p\}))\qtq{and}
\pi_1((\p^1\times P^0)\cap F^{-1}(U), (0{\colon\!} 1)\times \{p\})
$$
have the same image in $\pi_1(U,x_0)$. The latter image has finite index
by (\ref{maps.fund.grps}), thus $g:=F_p:\p^1\to X$
is as desired.

In order to finish, it is sufficient to prove
that if  $g_1, g_2:\p^1\to X$ are as in (\ref{biggest.sg.lem}.2)
then there is a third morphism
$g:\p^1\to X$ such that
$$
\im[\pi_1(g^{-1}(U),(0{\colon\!} 1))\to \pi_1(U,x_0)]\supset
\im[\pi_1(g_i^{-1}(U),(0{\colon\!} 1))\to \pi_1(U,x_0)]
$$
for $i=1,2$. 
To do this, let $S=\{1,2\}$ be a 2 point scheme
and take $f_0:\p^1_S\to X$ to be $g_i$ on $\{i\}\times \p^1$. 
Construct $h:Y\to D$ and $F:Y\to X$ as in (\ref{smoothing.lem}).
Set $W:=F^{-1}(U)$ and apply (\ref{fg.uppsem.lem}) twice 
for $i=1,2$ 
with $C_0:=f_i^{-1}(U)$.
For any $d\in D^0$,  $g:=F_d$ has the required property.\qed

\begin{say}[Proof of (\ref{main.thm})]\label{pf.main.thm} {\ }

Let $H<\pi_1(U,x_0)$ be the subgroup obtained in
(\ref{biggest.sg.lem}). We are done if $H=\pi_1(U,x_0)$.
Otherwise there is a corresponding irreducible \'etale cover
$(x'_0\in U')\to (x_0\in U)$. Let $w:X'\to X$ be the normalization
of $X$ in the function field of $U'$.
$X$   is simply connected by (\ref{rc->simp.conn}),
 hence by the purity of branch loci
(cf.\ \cite[X.3.1]{sga1})
there is a divisor $D'\subset X'$ such that $w$ ramifies along $D'$.
Let $D\subset X$ be the image of $D'$. 
By construction $D\subset X\setminus U$.
We derive a contradiction as follows.

Let $f:\p^1\to X$ be a morphism such that
$f(0{\colon\!} 1)=x_0$ and $f(\p^1)$ intersects
$D$ transversally at a point $x_1:=f(1{\colon\!} 1)$.
This implies that the local fundamental group of
 $\p^1\setminus\{(1{\colon\!} 1)\}$
at $(1{\colon\!} 1)$ 
surjects onto 
the local fundamental group of $X\setminus D$ at $x_1$.
(See \cite{g-murre} for the  local fundamental group 
of a divisor in a  variety. There it is called 
 the fundamental group
of the formal neighborhood of a divisor.)
Therefore, if 
 $f$ lifts to $f':\p^1\to X'$ then
 $X'\to X$ is \'etale at $f'(1{\colon\!} 1)$.
If $X'\to X$ is a Galois extension the we already have a contradiction
since $X'\to X$ ramifies everywhere above $D$.
In the non--Galois case 
$p$ may be unramified along some of the irreducible components of
$p^{-1}(D)$, and we only get a contradiction if the
image of 
 $f'$  intersects  $D'$.

If $X'\to X$ is not Galois, we need to proceed in a
somewhat roundabout way. First I give the outlines; a 
precise version is given 
afterwards.

It is clear from the definition that there are many
maps $\p^1\to X'$, thus $X'$ 
(or rather any desingularization of $X'$) 
 should be rationally connected.
 This indeed follows from  (\ref{lifting.lem})
applied to any $F:\p^1\times P\to X$  as in 
 (\ref{rc.equiv.def}.4).  Thus by
 (\ref{rc.equiv.def}.5) there is a morphism $f':\p^1\to X'$
 which passes through $x'_0$ and intersects $D'$
at a smooth point $x'_1$.
 We obtain a contradiction if $w\circ f'$ is the limit of a sequence of maps
  $f_t:\p^1\to X$ such that
  \begin{enumerate}
  \item $f_t$ lifts to $f'_t:\p^1\to X'$ and $f'$ is
 the limit of the maps $f'_t$,
  \item the image of $f_t$ intersects $D$ transversally.
  \end{enumerate}

First we prove a general  lifting property for families of morphisms
and then we proceed to construct the morphism $f'$.

\begin{lem}\label{lifting.lem}
 Let $V$ be a normal variety and $G:\p^1\times V\to X$ a morphism
such that $ G((0{\colon\!} 1)\times V)=\{x_0\}$. Assume that
$H^1(\p^1, G_v^*T_X(-2))=0$ for some $v\in V$. 
Then $G$ can be lifted to $G':\p^1\times V\to X'$.
\end{lem}

Proof.  First we show that such a lifting exists over an open subset
of  $\p^1\times  V$.
Choose $V^0\subset V$
 such that $H^1(\p^1, G_v^*T_X(-2))=0$ for every $v\in V^0$
(this is possible by the upper semi continuity of cohomology groups)
and such that we have an exact sequence
$$
\pi_1(G_v^{-1}(U), (0{\colon\!} 1)\times v)\to 
\pi_1((G^0)^{-1}(U),(0{\colon\!} 1)\times v)\to \pi_1(V^0,v)\to 1
$$
for every $v\in V^0$
where $G^0:=G|_{\p^1\times V^0}$  (this is possible by (\ref{open.fibr.lem})). 
Then
$$
 \im[\pi_1((G^0)^{-1}(U), (0{\colon\!} 1)\times v)\to \pi_1(U,x_0)]\subset H,
 $$
 thus $G^0|_{(G^0)^{-1}(U)}$ can be lifted to
 $G':(G^0)^{-1}(U)\to U'$. This extends to a rational map
  $G':\p^1\times V^0\map X'$.
  
  Next let $\Gamma\subset  (\p^1\times V)\times X$ be the graph of $G$
  and  $\Gamma^*\subset  (\p^1\times V)\times X'$  its preimage.
  The rational lifting $G'$ corresponds to an irreducible
  component $\Gamma'\subset \Gamma^*$ such that the projection
  $\Gamma'\to (\p^1\times V)$ is birational. Since $X'\to X$ is finite,
  so is $\Gamma^*\to \Gamma$. Thus $\Gamma'\to \Gamma\to \p^1\times V$
 is finite and
 birational,
  hence an isomorphism.\qed
  \medskip

Let $\phi:X''\to X'$ be any desingularization and $x'_1\in D'$
a point such that $D'$ and $X'$ are smooth at $x'_1$ and $\phi^{-1}$ is
a local isomorphism near $x'_1$. 
By (\ref{rc.equiv.def}.5) there is a dominant morphism
$F:\p^1\times P\to X'$ such that
\begin{enumerate}
\item $F((0{\colon\!} 1)\times P)=\{x'_0\}$, 
\item $F((1{\colon\!} 1)\times P)=\{x'_1\}$, and
\item the image of $F_p:\p^1\to X'$ is transversal
to $D'$ at $x'_1$ for every $p\in P$.
\end{enumerate}

Let us now consider the dominant morphism $w\circ F:\p^1\times P\to X$. 
By the first part of (\ref{free.maps.lem}), there is a $p_0\in P$ such that
$H^1(\p^1, (w\circ F_{p_0})^*T_X(-2))=0$. 
Thus by the second part of (\ref{free.maps.lem}) 
there is a pointed variety $q_0\in Q$ and a morphism
$G: \p^1\times Q\to X$ such that 
$G_{q_0}=w\circ F_{p_0}$ and 
$G$ is smooth away from
$(0{\colon\!} 1)\times Q$. In particular, $G^{-1}(D)\subset \p^1\times Q$
is a  generically  smooth divisor, hence  there is a dense open set
$Q^0\subset Q$ such that the projection
$G^{-1}(D)\to Q$ is smooth over $Q^0$.
This means that 
the image of $G_q:\p^1\to X$
intersects $D$ transversally for every $q\in Q^0$.
(Note that $G^{-1}(D)$ denotes the inverse image scheme.)

By (\ref{lifting.lem}), $G$ can be lifted to
$G':\p^1\times Q\to X'$. On  $\p^1\times \{q_0\}$ the lifting
agrees with  $F_{p_0}$,  hence $G'(\p^1\times  \{q_0\})$
intersects $D'$ at the point $x'_1$. This implies that
$(G')^{-1}(D')$ is a divisor and so by shrinking $Q$
we may assume that the image of $G'_q$ intersects $D'$ for every
$q\in Q$. Thus  $G_q=w\circ G'_q$ is never transversal to
$D$, a contradiction.

Condition (\ref{main.thm}.1) holds by construction and
a general choice of $f$ satisfies
(\ref{main.thm}.2) by \cite[II.3.14]{rcbook}.
\qed 
  \end{say}

\begin{say}[Proof of (\ref{main.cor})] {\ }

Pick $f_1:V_1\to U_{\bar K}$ defined over $\bar K$ such that
$\pi_1(V_1,0)\to \pi_1(U_{\bar K},x_0)$ is surjective
 and $H^1(\p^1, \bar f_1^*T_X(-2))=0$.
$f_1$ is defined over a finite Galois extension
$L\supset K$; let $f_i:V_i\to U_{\bar K}$ be its conjugates.
Each of these extends to a morphism $\bar f_i:\p^1\to X_{\bar K}$.
Let $S=\spec_KL$. We can view the morphisms $f_i$ as one morphism
$f_0:\p^1_S\to X$ defined over $K$. 
By (\ref{smoothing.lem}) we obtain  $h:Y\to D$
and $F:Y\to X$, all defined over $K$.
 Let $d\in D^0(\bar K)$ be any point. Then by 
(\ref{fg.uppsem.lem}) we see that
$$
\im[\pi_1(Y_d,0)\stackrel{F_d}{\to} \pi_1(U_{\bar K},x_0)]\supset
\im[\pi_1(V_1,0)\stackrel{f_1}{\to} \pi_1(U_{\bar K},x_0)],
$$
and the latter image is $\pi_1(U_{\bar K},x_0)$ by assumption.
$D^0$ is a Zariski open set in a curve $D$ with a smooth $K$-point
$0$. By (\ref{large.field.defn})
this implies that $D(K)$ is dense in $D$, hence $D^0(K)\neq \emptyset$.
By choosing $d\in D^0(K)$  we obtain an
 open
 subset $0\in V\subset \a^1$ and a morphism $f:V\to U$ 
 (all defined over $K$) such that $f(0)=x_0$ and
$$
\pi_1 (V_{\bar K},0)\onto \pi_1 (U_{\bar K},x_0)\qtq{is surjective.}
$$
The fundamental group of  a $K$-scheme $W$
is related to the fundamental group of   $W_{\bar K}$
by the exact sequence (cf. \cite[IX.6.1]{sga1})
$$
1\to \pi_1 (W_{\bar K},0)\to \pi_1 (W,0)\to \gal(\bar K/K)\to 1.
$$
This implies that
$$
\pi_1 (V,0)\onto \pi_1 (U,x_0)\qtq{is also surjective.}
$$
\qed
\end{say}

\begin{exmp} For every $n\geq 4$ here I give an example
 of a rational threefold $X$,
a normal crossing divisor $F\subset X$ and a smooth rational curve
$B\subset X$ such that $B$ intersects $F$ everywhere transversally,
$B$ intersects every irreducible component of $F$ and
the image of $\pi_1(B\setminus F)\to \pi_1(X\setminus F)$
has index $n$ in $\pi_1(X\setminus F)$.
 
Let us start with a similar surface example.
Let $g_1:\p^1\to \p^1$ be a degree $n$ morphism
with critical points $x_1,\dots,x_{2n-2}\in \p^1$ and different
 critical values $y_1,\dots,y_{2n-2}$,
 Choose 3 other points $x_{2n-1},x_{2n},x_{2n+1}\in \p^1$
 such that $g_1(x_{2n-1})$ and $g_1(x_{2n})$ are different critical
 values of $g_1$
 and $g_1(x_{2n+1})$ is not a critical value of $g_1$. 
 Let $g_2:\p^1\to \p^1$ be  a morphism with $g_2(x_i)=x_i$
 for  $i\leq 2n+1$.
 Consider the morphism $h:(g_2,g_1):\p^1\to \p^1\times \p^1$
 
Set $D=\p^1\times \{y_1,\dots,y_{2n-2}\}$ and
 $U=\p^1\times \p^1\setminus D$.

$\im h$ intersects every irreducible componnet of $D$ transversally 
at  $n-2$ points and the image of $\pi_1(h^{-1}(U))\to \pi_1(U)$
has index $n$ in $\pi_1(U)$. 

Here $\im h$ is also tangent to every irreducible component
of $D$. The tangencies can be resolved by 2n   blow ups,
but then the birational transform of $\im h$ does not intersect
every boundary component. To remedy this situation,
take another morphism
$g_3:\p^1\to \p^1$   with $g_3(x_i)=x_i$
 for  $i\leq 2n+1$.
Set $Y:=\p^1\times \p^1\times \p^1$, $F=D\times \p^1$ and
$H:=(g_3,g_2,g_1):\p^1\to Y$.
 Again the only problem is
that $\im H$ is tangent to every irreducible component
of $F$. These can be resolved by blowing up
 two suitable  smooth curves. 
First we take a smooth curve  
 which passes through
every point of tangency and also 
through the point $H(x_{2n-1})$.
After blow up, the birational transform of $\im H$
intersects every boundary component transversally but 
above each point there is a point common to two boundary components
and to the birational transform of $\im H$. 
Next take a smooth curve which passes through all these points
with a general tangent direction there and also through 
 $H(x_{2n})$. We can also assume that neither of the two curves
passes through $H(x_{2n+1})$. 
Doing two such blow ups creates two new boundary components and
the birational transform of $\im H$   intersects both of them.
The fundamental group computation is unchanged.

Varying $g_2,g_3$ we obtain many morphisms all of which pass through
the point $H(x_{2n+1})$. 
\end{exmp}

\begin{ack}  I   thank J.-L.\ Colliot-Th\'el\`ene and D.\ Madore
for many helpful
comments and references.
Partial financial support was provided by  the NSF under grant number 
DMS-9970855. 
\end{ack}

\vskip1cm

\noindent Princeton University, Princeton NJ 08544-1000

\begin{verbatim}kollar@math.princeton.edu\end{verbatim}

\end{document}